\documentclass[10pt]{article}

\usepackage[utf8x]{inputenc}
\usepackage[english]{babel}

\usepackage{amsmath}
\usepackage{physics}
\usepackage[margin=1.5in]{geometry} 
\usepackage{setspace}
\usepackage{mathtools}
\usepackage{relsize}
\usepackage{amsmath,amsthm,amssymb,amsfonts}
\usepackage{latexsym}
\usepackage{tikz-cd}
\usepackage[square,sort,comma,numbers]{natbib}
\usepackage{graphicx}
\usepackage{url}

\title{What is the Perfect Shuffle?}
\author{James Enouen}
\date{September 2019}

\begin{document}

\maketitle

\begin{abstract}
When shuffling a deck of cards, one probably wants to make sure it is thoroughly shuffled.
A way to do this is by sifting through the cards to ensure that no adjacent cards are the same number, because surely this is a poorly shuffled deck.
Unfortunately, human intuition for probability tends to lead us astray.
For a standard 52-card deck of playing cards, the event is actually extremely likely.
This report will attempt to elucidate how to answer this surprisingly difficult combinatorial question directly using rook polynomials.
\footnote{A special thanks to Srivatsa Srinivas for connecting Gessel's polynomial to this problem.}
\end{abstract}

\section{Introduction}
We will say that a shuffle of a standard 52-card deck is a \emph{perfect shuffle} if any pair of adjacent cards in the deck have a different value from one another.

\begin{figure}[h!]
\centering
\includegraphics[scale=.08]{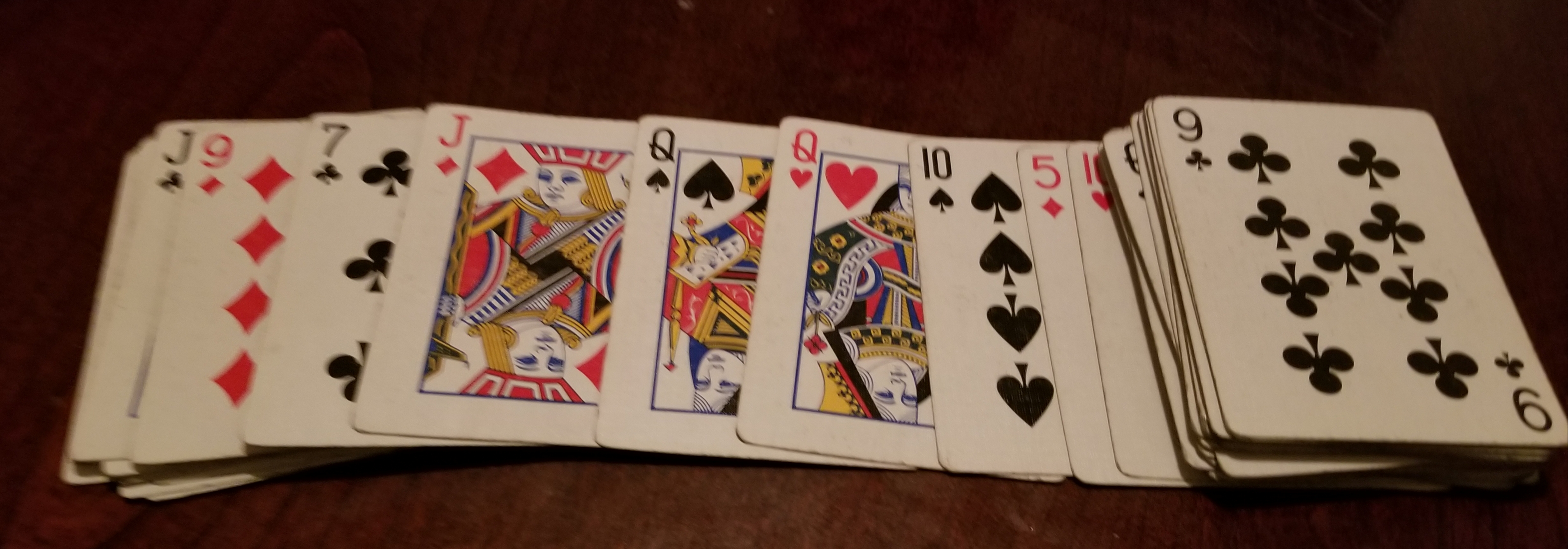}
\caption{An imperfect shuffle}
\label{fig:universe}
\end{figure}

Formally, we can see this as a permutation on 52 elements where the first four elements are of a first color, the second four elements are of a second color, and on until the last four elements are of the thirteenth color.
Answering this problem was looked at in a 2013 IJPAM article by Yutaka Nishiyama \citep{nishiyama13}.
Unfortunately, the perspective he took to analyze the problem became very computationally intensive for decks with 52 cards, so he used Monte Carlo methods to approximate the answer for the 52-card deck\citep{nishiyama13}.
The only other places I have found the exact answer to this question is on a swedish forum from 2009 \citep{swedish} and a french blog from 2014 \citep{french}.
As far as I can tell, the two who answered these questions computed the answer with brute force using more efficient code and more resources than Nishiyama's attempt.
Here, we address the problem using Ira Gessel's 1988 generalization of rook polynomials to achieve a solution which is much less computationally restrictive \citep{gessel88}.

\section{Introduction to Rook Polynomials}
Rook polynomials were developed studying the number of ways to place rooks on a chessboard.
Our study will fairly closely follow a combination of Gessels' work in \citep{gessel88} and \citep{gessel13}.
For a given size of chessboard, the rook polynomial counts the number of ways to place differing amounts of non-attacking rooks on that chessboard.
Let $n \in \mathbb{N}$ and let $[n]$ denote $\{1,...,n\}$.
We consider our chessboard to be $[n] \times [n]$.  Let a \emph{board}, B, be a subset of the chessboard $B \subseteq [n] \times [n]$.  
For this board, we define the rook number $r_k(B)$ to be the number of ways to put $k$ rooks onto this board such that none of them are `attacking' the other (none are in the same row or column).  We will take $B' = \{(2,2), (3,2), (3,3)\}$ for our illustrative examples:

\begin{figure}[h!]
\centering
\includegraphics[scale=.23]{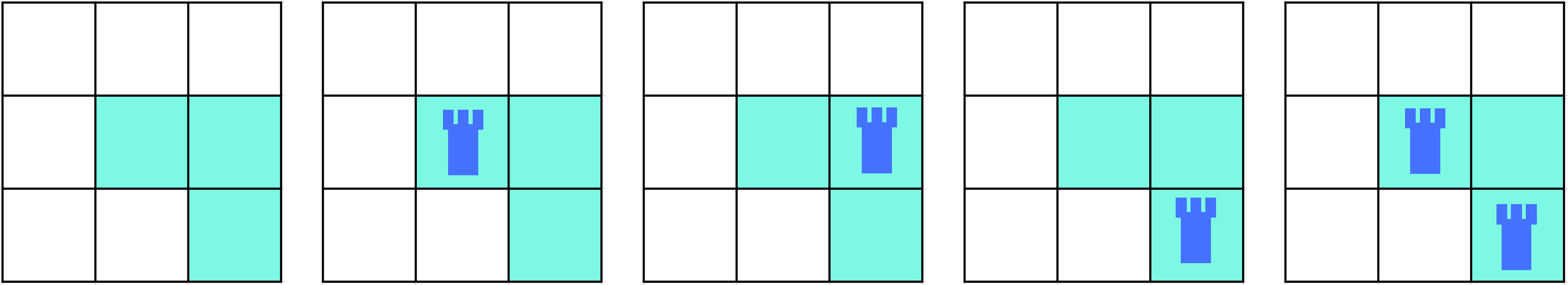}
\caption{We see: $r_0(B') = 1$, $r_1(B') = 3$, and $r_2(B') = 1$.}
\label{fig:rooks}
\end{figure}

Let $S_n$ denote the set of permutations of $[n]$.  
We can associate each $\pi \in S_n$ to another subset of our chessboard.
We do this with the set $\{ (i,\pi(i)) : i \in [n] \}$.
We can now see how many `hits' each permutation has with a given board.
Formally, we can define a `hit function'
\[ h_B : S_n \rightarrow \mathbb{N}_0 = {0,1,2,...}\]
\[ h_B(\pi) := | \{ (i,\pi(i)) : i \in [n] \} \cap B | \]
The associated hit numbers count how many of the $n!$ permutations hit that many times.
\[ h_k(B) := | \{ \pi \in S_n : h_B(\pi) = k \} | \hspace{17pt} \textrm{for} \hspace{3pt} k \in \mathbb{N}_0\]

\begin{figure}[h!]
\centering
\includegraphics[scale=.42]{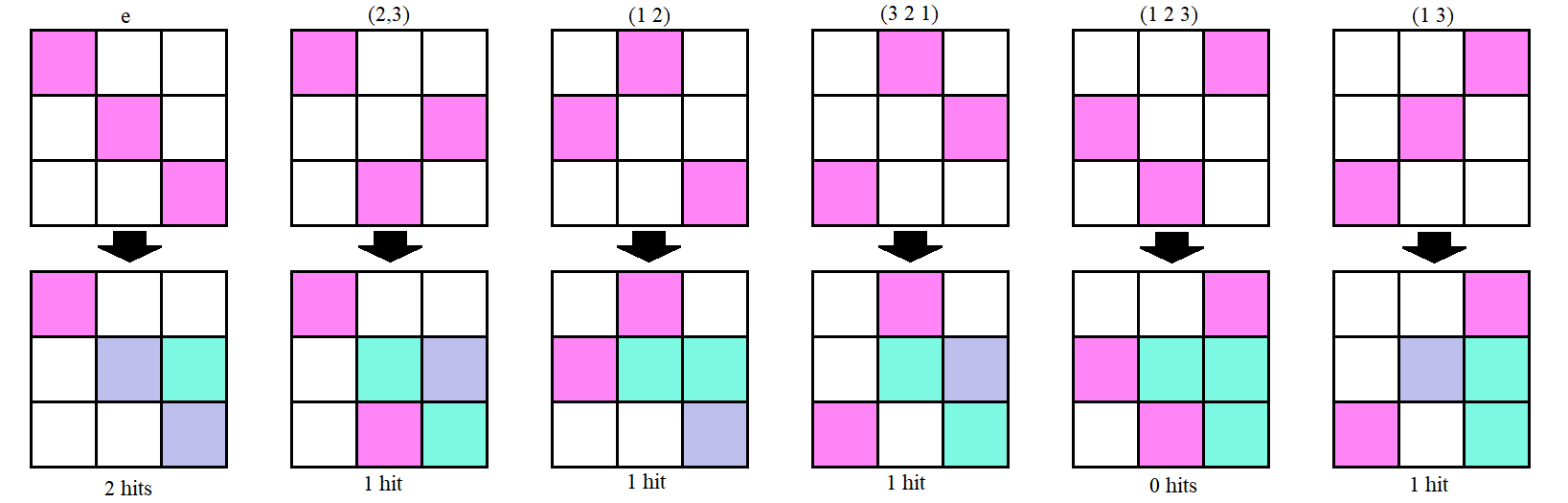}
\caption{We see: $h_0(B') = 1$, $h_1(B') = 4$, $h_2(B') = 1$.}
\label{fig:hits}
\end{figure}

\newpage

We will now see the identity which relates the rook numbers and the hit numbers by considering arrangements of rooks as partial permutations of $[n]$.
\[ \mathlarger{{\sum_i}} h_i(B) \binom{i}{j} = r_j(B)\cdot(n-j)! \hspace{17pt} \forall j\in \mathbb{N}_0 \]
\begin{proof}
This equality will be achieved by double counting the number of pairs $(\pi, H)$ where $\pi$ is a permutation and $H$ is a j-subset of the set of $\pi$'s hits ($\{ (i,\pi(i)) : i \in [n] \} \cap B$)

The left-hand side picks $\pi$ first.  
Call i the number of hits of $\pi$, $i=h_B(\pi)$, and then take all $\binom{i}{j}$ subsets of the i hits which have size j.
Since there are $h_i(B)$ permutations which have hit number i and contribute $\binom{i}{j}$ pairs to our sum, we have ${\sum_i} h_i(B) \binom{i}{j}$ pairs in total.

\begin{figure}[h!]
\centering
\includegraphics[scale=.36]{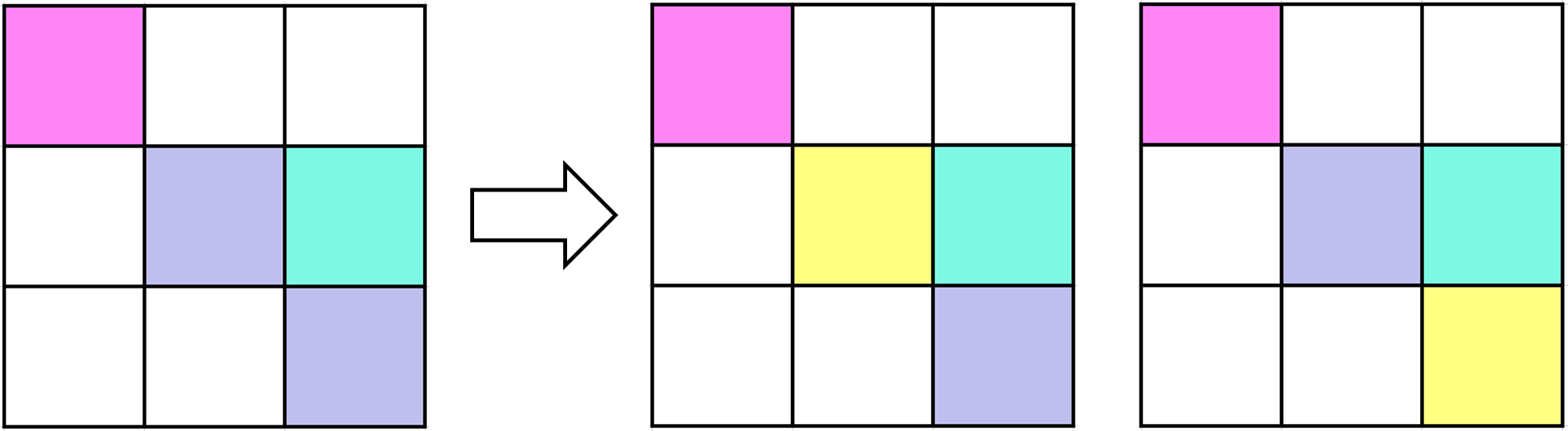}
\caption{We first pick $\pi = e$ and then choose the $\binom{2}{1}$ different subsets (of size 1) of its 2 hits.}
\label{fig:LHS}
\end{figure}

The right-hand side picks H first.
We know that there are $r_j(B)$ subsets of the board which have size j.
For each of these, we can then extend them to a permutation by filling the empty rows and columns in (n-j)! ways.

\begin{figure}[h!]
\centering
\includegraphics[scale=.36]{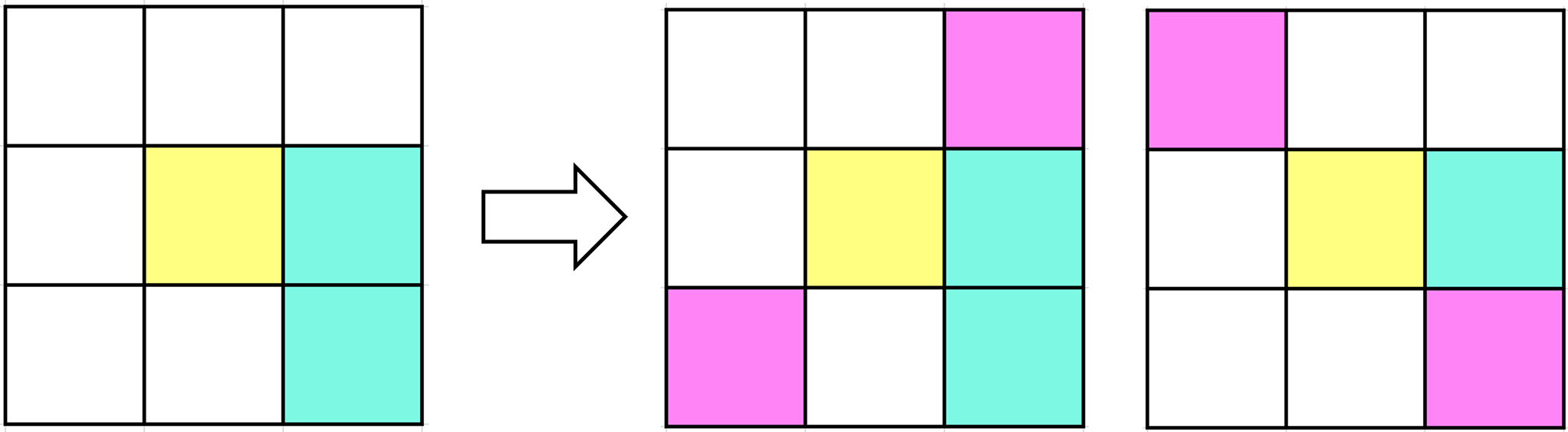}
\caption{We first pick $H=\{(2,2)\}$ and extend to the $(3-1)!$ different permutations}
\label{fig:RHS}
\end{figure}

\end{proof}

\subsection{Defining the Rook Polynomial}

Given this identity for all j, we can multiply each of these identities by $t^j$, yielding:
\[ \mathlarger{{\sum_i}} h_i(B) \binom{i}{j} t^j = r_j(B)\cdot(n-j)!\cdot t^j \hspace{17pt} \forall j\in\mathbb{N}_0\]
Summing over j yields
\[ \mathlarger{{\sum_i}} h_i(B) \mathlarger{{\sum_j}} \binom{i}{j} t^j = \mathlarger{{\sum_j}} r_j(B)\cdot(n-j)!\cdot t^j \]
\[ \mathlarger{{\sum_i}} h_i(B) (1+t)^i = \mathlarger{{\sum_j}} r_j(B)\cdot(n-j)!\cdot t^j\]
Plugging in $t=-1$ yields
\[ h_0(B)  = \mathlarger{{\sum_j}} (-1)^j\cdot r_j(B)\cdot(n-j)! \]

This number $h_0(B)$ counts how many permutations totally avoid our subset B.
This is the critical identity which fuels the study of rook polynomials, but this identity can equally be derived through the principle of inclusion and exclusion.

Regardless, corresponding to this equation, let us define the rook polynomial:
\[ r_B(x) := \mathlarger{{\sum_k}} (-1)^k\cdot r_k(B)\cdot x^{n-k} \]
Let $\phi$ be a linear functional on polynomials in x with the effect: 
\[ \phi(x^k) = k! \]
Thus, $h_0(B) = \phi(r_B(x))$ which may seem a little convoluted at first, but it results in the following fantastic property:

$r_{B_1}(x) * r_{B_2}(x) = r_{B_1\oplus B_2}(x).\hspace{17pt}$ where $B_1\oplus B_2$ is the direct sum of two boards as depicted below.

\begin{figure}[h!]
\centering
\includegraphics[scale=.36]{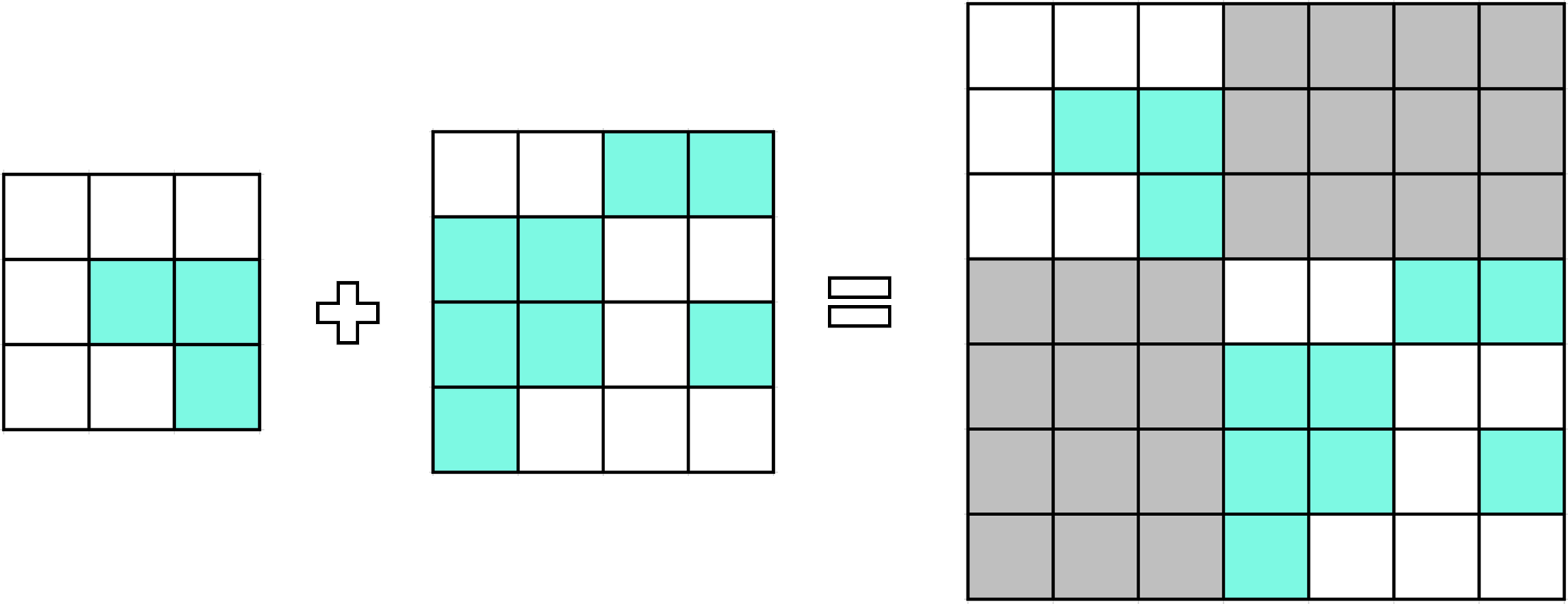}
\caption{Direct sum of two boards}
\label{fig:comp}
\end{figure}

\subsection{Implications of the Product Formula}
Take $l_n(x)$ to be for the complete board
$l_n(x) = r_{[n]\times[n]}(x)$.
We need to count the number of ways to put k rooks on the nxn board.
The first rook has $n^2$ places to go, the second will then have $(n-1)^2$ places, then $(n-2)^2$, etc.  Since we picked these with respect to order, we need to divide by the $k!$ different orders to get the number of ways and then:
\[
\frac{(n)^2\cdot(n-1)^2\cdot ... \cdot(n-(k-1))^2 }{k!} = k! \cdot \frac{((n)\cdot(n-1)\cdot ... \cdot(n-(k-1)))^2}{(k!)^2} = k! \cdot\binom{n}{k}^2 
\]
So, $l_n(x) = \sum_{k=0}^n (-1)^k \binom{n}{k}^2 k! x^{n-k}$.  This allows to write the solution for the number of `generalized derrangements.'

The number of permutations of $n = n_1 + ... + n_r$ objects where $n_i$ objects have the color i such that i and $\pi(i)$ have different colors is:
\[
\phi\bigg{(}\mathlarger{\prod_{i=1}^r} l_{n_i}(x) \bigg{)}
\]
This fact simply follows from our product of boards identity by using the full boards $[n_i]\times[n_i]$ and full rook polynomials $l_{n_i}(x)$.
(The case of derrangements is $n_i = 1$ for all $i\in[r]$.)
This result was proved by Evens and Gillis in 1976, without using the connection to rook theory.

\begin{figure}[h!]
\centering
\includegraphics[scale=.33]{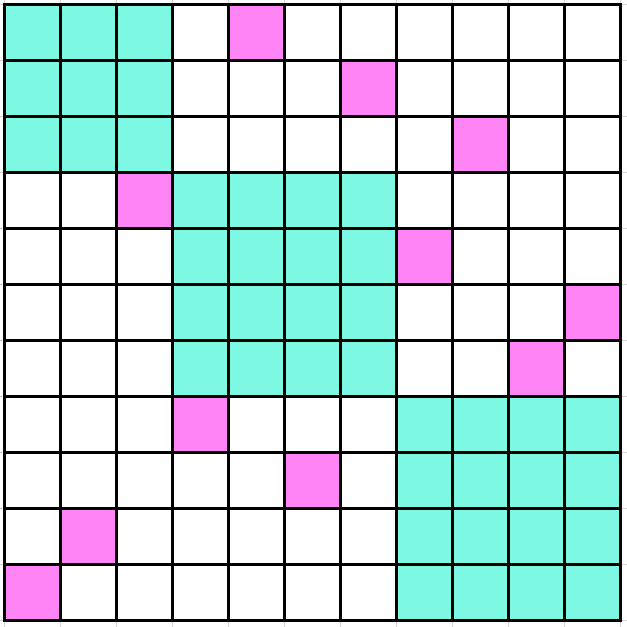}
\caption{An example of a generalized derrangement with $n_1=3,n_2=4,n_3=4$}
\label{fig:der}
\end{figure}

\section{Generalized Rook Polynomials}
To generalize this rook polynomial beyond permutations of $[n]$ and the conditions $\pi(i) = j$, we can use Ira Gessel's \citep{gessel88} definition of a generalized rook polynomial.

Let us have sets $T_0,T_1,T_2,...$ which have cardinalities $M_0,M_1,M_2,...$ .
For each $n\in\mathbb{N}_0$ we will additionally have a ``set of conditions'' $C_n$ which satisfy $C_0 \subseteq C_1 \subseteq C_2 \subseteq ...$.  
To each condition $c\in C_n$ we will have a set $T_n^c \subseteq T_n$ which are the elements of $T_n$ which satisfy the condition.
For a set of multiple conditions $A\subseteq C_n$ we will say the set satisfying all of these properties is $T_n^A = \cap_{a\in A} T_n^a$.

The following property is what we need to stay in the rook polynomial structure:
If $A\subseteq C_n$, one of the following occurs:
\vspace{-5mm}
\begin{itemize}
  \item $T_m^A = \emptyset$ for all $m\geq n$. ``A is incompatible''
  \item There is a $\rho(A)\in\mathbb{N}_0$ such that for every $m\geq n$ there is a bijection $T_m^A\rightarrow T_{m-\rho(A)}$
\end{itemize}

Another technical condition we need is that our sequence $(M_n)\in\mathbb{N}_0$ does not satisfy any linear homogeneous recurrence equation.
This condition is needed so that $\rho(A)$ will be uniquely determined when it exists.

For a set of conditions in $C_n$, take $B\subseteq C_n$.
This means $B\subseteq C_m$ for all $m\geq n$.
We want to count the elements of $T_m$ for $m\geq n$ which satisfy \textit{none} of the conditions in B.
We will denote this set $T_m/B$.
By inclusion-exclusion, we have that
\[ T_m/B = \sum_{\substack{A\subseteq B \\ compatible}} (-1)^{|A|}\cdot M_{m-\rho(A)}    \]
Correspondingly, we will finally define our ``generalized rook polynomial'' to be:
\[ r_B(x) = \sum_{\substack{A\subseteq B \\ compatible}} (-1)^{|A|}\cdot x^{n-\rho(A)}    \]

If we define $\phi(x^n) = |T_n| = M_n$, then for all $m\geq n$, $\hspace{17pt} |T_m/B| = \phi( r_B(x) \cdot x^{m-n} )$.
Because of our linear recurrence restriction on $(M_n)$, this equation is actually able to uniquely determine $r_B(x)$.

This new definition still has the property that the product of the polynomials for two disjoint `boards' are equal to the the polynomial of their disjoint attachment.
This will again be useful to us.

\subsection{Connection to the Original Rook Polynomial}

We will first tie this more general definition back to our original setting.
First, our sets $T_n$ were the sets of permutations of $[n]$, so $T_n = S_n$ for all $n\in\mathbb{N}_0$, hence $M_n =n!$ for all $n\in\mathbb{N}_0$.
The conditions are slightly more tricky.
These conditions are associated to the `boards' which were subsets of $[n] \times [n]$ and prescribed which spots on the board the permutations were not allowed to `hit.'
So, $C_n = \{ ``\pi(i) = j''$ for $(i,j)\in [n]\times[n] \}$ and a set of conditions in $C_n$ $B\subseteq C_n$ is the same as a board $B\subseteq [n]\times [n]$

We can now check if our original setting indeed satisfies the property of compatibility claimed to be essential.
Let $A \subseteq C_n$.
Suppose there is a distinct pair $(i,j),(i',j')\in A$ with $i=i'$ or $j=j'$, then $T_m^a = \emptyset$ and A is incompatible, because there is no permutation which sends one element to two different values and there is no permutation which sends two elements to the same value.
Otherwise, there is a pretty simple bijection by seeing that each condition in A fixes exactly one input and output of our permutation.
Each time we fix an element of our permutation of n elements, we are left with a permutation of (n-1) elements.
This is the bijection $S_m^A \rightarrow S_{m-|A|}$ for all A which are compatible.
This means that $\rho(A) = |A|$.

Now we can see that our definitions now align because the compatible subsets of size k correspond exactly to placements of non-attacking rooks onto the board:
\[ \mathlarger{{\sum_k}} (-1)^k\cdot r_k(B)\cdot x^{n-k} = r_B(x) = \sum_{\substack{A\subseteq B \\ compatible}} (-1)^{|A|}\cdot x^{n-\rho(A)}    \]

\subsection{Linear Permutations}

In this setting, we will also take $T_n = S_n$ and $M_n = n!$, but our conditions will be very different.
Our elementary conditions will be ``i is immediately followed by j'' where we see a permutation as the one-line notation.
Equivalently, this condition is ``for some k, $\pi(k) = i$ and $\pi(k+1) = j$.''
So our set of conditions are $C_n = \{$``i is immediately followed by j'' $: i\neq j\in [n] \} $ and $C_0\subseteq C_1\subseteq C_2\subseteq ...$ is $\{\}\subseteq \{\}\subseteq \{$``1 follows 2'', ``2 follows 1'' $\}\subseteq ... $

We will now show that the compatibility property holds.
Let $A\subseteq C_n$
Similar to before, if we have two distinct conditions $(i,j),(i',j')\in A$ with $i=i'$ or $j=j'$, then A must be incompatible because $i=i'$ can't be immediately followed by two different numbers $j,j'$ or $j=j'$ can't be immediately preceded by two different numbers $i,i'$.
Otherwise, the bijection to a permutation of a smaller number of elements comes by viewing each adjacent pair i,j as a single element.
Each time we group two elements of our permutation of n elements, we are left with a permutation of (n-1) elements (where one element consists of two).
For larger strings of pairs i,j,k,l,... which are all adjacent we can see these r elements all as one element or we can individually put together pairs of elements.
Regardless, we ultimately see that we have a bijection $S_m^A \rightarrow S_{m-|A|}$ for all A which are compatible.
This means that $\rho(A) = |A|$.

\subsection{Product Formula Application}

We will again take the full set of conditions in order to yield a useful formula.
Take $l_n^*(x) = r_{C_n}(x)$ where $C_n$ is defined as above.

We only need to look at compatible subsets, so we only need to pick elements in different rows and columns if we see $C_n$ as the chessboard without the diagonal.
We will now choose a compatible set of conditions of size k.
The first condition we choose has $n^2-n = n(n-1)$ places to go, the second will then have $(n-1)^2-(n-1) = (n-1)(n-2)$ places, then $(n-2)(n-3)$, etc.
Since we picked these in a particular order, we need to divide by the $k!$ to correctly count the number of subsets of size k.

\[
\frac{1}{k!}\cdot(n)(n-1)\cdot(n-1)(n-2)\cdot ... \cdot(n-(k-1))(n-k) = k! \cdot \frac{(n)\cdot ... \cdot(n-(k-1))}{k!} \cdot \frac{(n-1)\cdot ... \cdot(n-(k))}{k!}\]
\[
= k! \cdot\binom{n}{k}\cdot\binom{n-1}{k}
\]

So, $l_n^*(x) = \sum_{k=0}^n (-1)^k \binom{n}{k}\cdot\binom{n-1}{k} k! x^{n-k}$.  
This allows to write the solution for the number of a specific type of linear permutation.

The number of linear arrangements of $n = n_1 + ... + n_r$ objects where $n_i$ are of color i, such that every adjacent object has a different colors is:
\[
\phi\bigg{(}\mathlarger{\prod_{i=1}^r} l_{n_i}^*(x) \bigg{)}
\]

This work on generalized rook polynomials has many more applications.
There are weighted sums instead of simple counting and the polynomials are used on other sets than just permutations.
If you are interested in this style of combinatorics, I highly encourage you to check out Gessel's paper \citep{gessel88}.
If you are not so interested, however, you are in luck because the linear permutation case is the solution to our original problem.

\section{Solution}
Our original problem asked about a set of 52 objects/ cards.
We asked that our permutation had no adjacent object of the same value.
We can see that we have 13 different values which we can now see as colors to see that we have 13 colors each with 4 associated objects/ cards.
So, we can take $r = 13$ and $n_i = 4$ for each $i \in [r]$.  
We first want to calculate
\[
\mathlarger{\prod_{i=1}^r} l_{n_i}^*(x) = (l_4^*(x))^{13}
\]
\begin{multline}
(l_4^*(x))^{13} = \\
-876488338465357824 x^{13} + 17091522600074477568 x^{14} - 159520877600695123968 x^{15} \\
+ 949054268820802240512 x^{16} - 4044281535242623254528 x^{17} + 13151570567369808936960 x^{18} \\
  - 33954920849889627734016 x^{19} + 71502295779064701517824 x^{20} - 125212657768448227540992 x^{21} \\
 + 185006084370341623234560 x^{22} - 233228682051017005596672 x^{23} + 253073982060156904538112 x^{24} \\
- 238025750670961148952576 x^{25} + 195147037097635696607232 x^{26} - 140102373840493649854464 x^{27} \\
+ 88405409991914856382464 x^{28} - 49175456453520166748160 x^{29} + 24169421980306186960896 x^{30} \\
- 10514786687648809353216 x^{31} + 4054104097647470051328 x^{32} - 1386375667685767249920 x^{33} \\
+ 420612294417061773312 x^{34} - 113190888701156917248 x^{35} + 27000049659200077824 x^{36} \\
- 5701677221962874880 x^{37} + 1063971192922619904 x^{38} - 175008802134196224 x^{39} \\
+ 25291193280417792 x^{40} - 3197671558907904 x^{41} + 351835440473088 x^{42} \\
- 33462483664896 x^{43} + 2727515172096 x^{44} - 188444475648 x^{45} \\
+ 10878057216 x^{46} - 514605312 x^{47} + 19420128 x^{48} \\
- 561912 x^{49} + 11700 x^{50} - 156 x^{51} + x^{52} \\
\end{multline}

So,
\[ \phi(\hspace{1pt} (l_4^*(x))^{13} ) =\]\[ 3,668,033,946,384,704,437,729,512,814,619,767,610,579,526,911,188,666,362,431,432,294,400 \]
Hence, the probability of a perfect shuffle is:
\[ \frac{\phi(\hspace{1pt} (l_4^*(x))^{13} ) }{52!} = 
\frac{672,058,204,939,482,014,438,623,912,695,190,927,357}{14,778,213,400,262,135,041,705,388,361,938,994,140,625} \approx 0.045476282331.\]
This means the chance of two adjacent cards being the same value is about $95.45\%$.
Interestingly, in the above probability, the numerator is prime and the denominator is $3^{5}\cdot 5^{10}\cdot 7^7\cdot 11^3\cdot 13^3\cdot 17^3\cdot 19^2\cdot 23^2\cdot29\cdot31\cdot37\cdot41\cdot43\cdot47$ which is always $p^{\lfloor{\frac{51}{p}}\rfloor}$ except for the lower prime factors of 2,3 which points to some small degree of symmetry in the space of ``perfect shuffles''.

\section{Conclusion}
The chance of a perfect shuffle is $\approx            4.5476282331\%$ and the chance of not is $\approx 95.4523717669\%$.
This easy to formulate question had a surprisingly sophisticated but rather elegant solution.
There are some questions which are obvious future directions of this problem.
The first is to consider how our probability changes when we consider the first and last cards of the deck to be `adjacent' to one another so that our deck of cards becomes a cyclic object.
The second is to consider instead of only $h_0(B)$ for our situation where we count the number of permutations satisfying none of the properties, we instead count each number which satisfy k of the properties (corresponding to $h_k(B)$.)
This will then give a distribution over the $52!$ permutations which count how many pairs of adjacent cards have the same value for a given shuffle.
Additionally, both of these questions can be asked simultaneously to give a distribution over the cyclic shuffles.
Hopefully this exposition was sufficient to understand the proof behind the coveted `probability of a perfect shuffle' and hopefully these future questions find their own answers as well.

\setstretch{0.8}
\bibliographystyle{plain}
\bibliography{references}

\end{document}